\documentclass{article}
\usepackage{amsmath}
\usepackage{theorem}
\usepackage{amssymb}
\usepackage{verbatim}
{\theorembodyfont{\rmfamily} \newtheorem{theorem}{Theorem}[section]
\newtheorem{proposition}{Proposition}[section]

}
{ \theorembodyfont{\rmfamily} }
\numberwithin{equation}{section}
\begin{document}
\title{DeBruijn Strings, Double Helices, and the Ehrenfeucht-Mycielski Mechanism}
\author{Terry R. McConnell\\
Syracuse University}
\maketitle
\begin{abstract}
We revisit the pseudo-random sequence introduced by Ehrenfeucht and Mycielski and its connections with DeBruijn strings.
\end{abstract}

\section{Introduction}

DeBruijn strings of order $n$ are finite strings of binary digits of minimal length that contain as substrings all $2^n$ possible binary strings of length $n$. Since the addition of a single binary digit can produce at most one new substring of length $n$, it is easy to see that a DeBruijn string must have length at least $2^n + n - 1.$ It is remarkable that this minimum is achieved - it is possible to pack the greatest possible variety of binary strings of given length into the minimum possible space.

For example, the string 0011101000 is a DeBruijn string of order 3, and the string 1111011001010000111 is a DeBruijn string of order 4. It follows from classical work of DeBruijn that there are exactly $2^{2^{n-1}}$ distinct DeBruijn strings of order $n$. (See, e.g., p. 136 of \cite{golomb}.) The would-be DeBruijn-string enthusiast can do no better than to begin with the excellent survey \cite{fred} of Fredrickson.

Here is a simple procedure to generate an example of a DeBruijn string of any desired order $n$: begin with an n-string that is all zeros except for a final 1. At each subsequent stage the next digit is a 1 unless that would cause a repeat of the terminal n-string. Otherwise, it is a zero. Continuing in this way will always produce a DeBruijn string that ends with $n$ zeros. (For a proof, see, e.g., pages 4-6 of \cite{me}, where the algorithm was naively presented as new. In fact, as noted in \cite{fred}, the algorithm is very well-known, often rediscovered, and goes back at least to 1934. It is known in the literature as the `prefer-one algorithm'.)

The Ehrenfeucht-Mycielski  (hereafter EM) sequence is an infinite binary sequence with pseudo-random properties first introduced in \cite{EM}. The sequence begins 010, and thereafter each additional digit is generated by an algorithm that attempts strenuously to avoid repetition. The precise algorithm is as follows: Find the longest suffix of the sequence generated so far that occurred  at least once previously. The next digit is then the opposite of the one that followed the penultimate occurrence of that suffix. Thus, the first 15 digits of EM are 010011010111000. It is easy to program a computer to calculate long prefixes of EM. For example, the first million digits are available at \cite{data}.

Since its introduction in 1992, the properties of EM have been studied by several authors. See, e.g., \cite{me},\cite{sutner},\cite{kieffer}, and \cite{herman}. It is sequence number A007061 in the Sloane On-line Sequence Encyclopedia \cite{sloan}.

The EM and prefer-one algorithms are similar in that both seek to avoid repetition. The main goal of this paper is to explore connections between the EM sequence and DeBruijn strings. For example, we show in section 3 below that if the comparison of binary strings in the EM algorithm is done in a finite buffer of fixed length $n$, then the resulting sequence will eventually become periodic, and that the repeating unit is a DeBruijn string of order $n + 1.$ This leads to a simple-minded and efficient algorithm for generating all possible DeBruijn strings of a given order. The algorithm may be new, albeit rather closely related to ALGORITHM 1 of \cite{fred}.

We also consider seeded variations of EM, and show, in particular, that the sequence uniquely determines the seed.

DeBruijn strings have many connections to other fields, including graph theory, coding theory, and shift-register sequences. For example, a DeBruijn string may be viewed as a Hamiltonian circuit of the DeBruijn graph of the same order. (See the following section for the relevant definitions.) If the edges of such a circuit are removed from the graph, the remaining graph always has a least 3 connected components, including simple loops on the vertices corresponding to the all-zero and the all-one string. In cases where there is only a single remaining (long loop) component, we term the original DeBruijn string a {\it double helix.} In the last section of the paper we show that double helices of all orders exist, and consider what happens when a double helix is used as seed for the EM sequence.

\section{Terminology and Notation}

The subject matter of this paper lies at the intersection of several fields, including Computer Science and parts of Mathematics that study combinatorics of words and free monoids. Accordingly, there is a variety of notation and terminology in common use. For example, we prefer the word `string' to `binary sequence' or `binary word', and will use it in the rest of the paper. Similarly, we prefer to use `tail' in place of `suffix', and `head' in place of `prefix'. We consider only strings over the alphabet $\{0,1\}$. It is likely some results can be generalized to larger alphabets, but we have not pursued this.

We shall use either small Greek letters for strings, or else capital Roman ones. The distinction is that Greek letters will be for `short' strings, and the others for `long' ones. String concatenation is indicated by juxtaposition of symbols. Individual binary digits will be denoted by $x,y,z,$ etc. A prime indicates complementary binary digit, i.e., if $x = 0$, then $x^{\prime} = 1$, and {\it vice versa.} The set of all binary strings is denoted $\Sigma.$

We denote by $\sigma_-$ the head, and by $\ _-\sigma$ the tail, of the string $\sigma$. More precisely, if $\sigma = x\tau y$ for some string $\tau$ and digits $x$ and $y$, then $\sigma_- = x\tau$, and $\ _-\sigma = \tau y.$  The length of $\sigma$ is $|\sigma|$. If $\sigma$ is a string of length at least $n$, then $\sigma_n$ shall denote the head of $\sigma$ of length $n$.  The symbols $\bar 0_n$ and $\bar 1_n$ indicate strings of length $n$ comprised of all 0s and all 1s respectively.

It will be convenient to use a few terms from \cite{me} in connection with the EM sequence. Given an initial string of EM, the longest tail $\sigma$ having an earlier occurrence will be called the {\it match string}, and $|\sigma|$ will be called the {\it match length.}

The DeBruijn graph of order $n$ is the directed graph ${\cal B}^n$ whose vertices are labeled by each of the possible binary strings of length $n$. If $\sigma = x\tau$ labels a given vertex, then there is an edge connecting from that vertex to each of the vertices labeled $\tau0$ and $\tau1$. There are no other edges than these. A given non-empty binary string X of length at least $n$ induces a path on ${\cal B}^n$ by starting at the vertex labeled with the initial n-string of X, and following with the vertices labeled by each successive substring of X of length $n$.

The graph ${\cal B}^{n+1}$ is, in a natural way, the dual graph of ${\cal B}^n$: Each edge $x\tau \to \tau y$ of ${\cal B}^n$ is made to correspond with vertex $x\tau y$ of ${\cal B}^{n+1}$. There is an edge in ${\cal B}^{n+1}$ between vertices corresponding to two given edges of ${\cal B}^n$, if and only if the latter edges meet at a common vertex. Under this correspondence, a given DeBruijn string X of order $n+1$ induces (i) a Hamiltonian path in ${\cal B}^{n+1}$; and (ii), an Eulerian path in ${\cal B}^n$. (In the first case, every vertex is visited exactly once; and, in the second case, every edge is crossed exactly once.)

\section{EM-variants and DeBruijn Sequences}

%\begin{lemma}
%An initial string of EM cannot be a square, i.e, an initial string cannot have the form %$\tau\tau$ for any nonempty string $\tau.$ In particular, an initial string cannot overlap with a later recurrence of itself.
%\end{lemma}

%Proof: Suppose an initial string of EM had the form $\tau\tau$. We may assume such $\tau$ is of %minimum length, necessarily greater than 1. Write $\tau\tau = \tau_{-}x\tau_{-}x$ for some %binary digit $x$. Since $\tau_{-}$ cannot overlap $\tau$ (by the minimality in the choice of %$\tau$,) the two copies of $\tau_{-}$ form a match. But this is clearly impossible. %$\hfill\square$

In this section we study certain variants of the EM sequence. First, let us modify the notion of string equality used in the definition of the EM sequence by interposing a transformation. Consider a function $h: \Sigma \to \Sigma$ which maps the empty string to the empty string. We shall deem two strings $\tau$ and $\sigma$ to be {\it h-equivalent} if $h(\tau) = h(\sigma).$ The EM sequence based on $h$, denoted $EM_h$, is the binary sequence $x_1x_2\dots$ defined by $x_0 =0$, and then as follows to define $x_{n+1}$ for $n \ge 0:$ Let $\sigma$ be the longest\footnote{This tail may be empty, in which case the latest previous occurrence is understood to be just before $x_n.$} tail of $x_1\dots{x_n}$ for which there is at least one $\tau \in x_1\dots{x_{n-1}}$ that is h-equivalent to $\sigma.$ Then $x_{n+1} = y^{\prime},$ where $y$ is the follower of the latest such $\tau.$

We also consider sequences with a given nonempty {\it seed string}, Z. One defines the binary sequence $EM_h(Z)$ by pre-pending Z to the strings $x_1\dots{x_n}$ and $x_1\dots{x_{n-1}}$ used in the definition of $EM_h$ above. If Z is nonempty then $x_1$ is determined by the same rule as is used for the later terms. (Note that $EM_h(Z)$ does not include the seed.)

If the function $h$ is one-one then $EM_h$ coincides with the usual Ehrenfeucht-Mycielski sequence, so we confine our attention to functions that are not one-one. A natural class of examples to consider are ones that model testing for string equality in a buffer of fixed finite length. This is, of course, the case in practise when a digital computer with finite memory is used to generate an initial segment of the EM sequence. Thus, given a fixed $n \ge 1$, define the sequence $EM_n(Z)$ ($EM_n$ if the seed is empty) by using for $h$ the function

\begin{equation}
h(\sigma) =
 \begin{cases}
 \sigma,& |\sigma| \le n\\
 \sigma_n,& \text{ otherwise.}
 \end{cases}
\end{equation}

The sequence $EM_n(Z)$ is identical to $EM(Z)$ until the match length first exceeds $n$. Computer experiments show that this initial section is followed by a chaotic stretch of seemingly indeterminate length, after which the sequence becomes periodic. The repeating section is a DeBruijn string of order $n+1$, as we show below.

\begin{proposition} For $n \ge 2$ a DeBruijn string X of order $n$ contains every string of length $n-1$ with opposite followers. All but the initial and terminal strings (which must be equal) occur exactly twice. The initial string occurs 3 times.
\end{proposition}

 Proof: Let $\tau$ be a string of length $n-1.$ Since both $\tau{0}$ and $\tau{1}$ occur in X, $\tau$ must occur at least twice. Suppose there where a third occurrence. Then $\tau$ must be the initial string of X to avoid a repetition of either $0\tau$ or $1\tau.$  But $\tau$ must also be the terminal string of X in order to avoid a repetition of either $\tau{0}$ or $\tau{1}.$\hfill$\square$

If a DeBruijn string Z of order $n$ begins with an n-string $\sigma$, then it ends with $\sigma_{-}.$ This follows from Proposition 3.1. Closely related to DeBruijn strings are {\it DeBruijn cycles.} Drop the final $\sigma_{-}$ of Z and view the resulting string as a cycle of length $2^n$ that begins afresh with the beginning $\sigma$ after $2^n$ terms. The literature on DeBruijn cycles is somewhat confusing since they have been called various things. For example, such cycles are called {\it full cycles} in \cite{fred}.

If each of the $2^n$ cyclic permutations of a DeBruijn cycle is extended to a DeBruijn string in the obvious way, then $2^n$ distinct DeBruijn strings are obtained, one for each possible initial $\sigma.$ It was originally proved by DeBruijn (see, e.g., p. 136 of \cite{golomb}) that there are exactly $2^{2^{n-1} - n}$ distinct DeBruijn cycles. Thus, there are exactly $2^{2^{n-1}}$ possible DeBruijn strings for a given $n$. (For $n = 4$, strings yielding 8 distinct DeBruijn cycles are listed on p. 135 of \cite{golomb}. These, together with their bitwise complements, comprise all 16 DeBruijn cycles of order 4. )

 For any given seed string Z, the sequence $EM_n(Z)$ can be generated by a finite state machine, and therefore it must eventually become periodic. Thus, unlike $EM(Z)$, the $EM_n(Z)$ are not transitive.\footnote{An infinite string is called {\it transitive} if it contains every possible string of finite length as a substring.}  On the other hand, they do contain all possible strings of length $n+1$:

 \begin{proposition} Every string of length $n+1$ occurs infinitely many times in $EM_n(Z).$
 \end{proposition}

 Proof: It is easy to see that $EM_n(Z)$ cannot eventually become constant, and therefore both binary digits occur infinitely many times. On the other hand, since $EM_n(Z)$ is ultimately periodic, there are sufficiently long strings that do not occur at all. Thus there is a shortest length string, $\tau$ say, that occurs at most finitely many times. Clearly $|\tau| \ge 2.$ Thus $\tau_{-}$ is a nonempty string occurring infinitely many times. Suppose $|\tau_{-}| \le n.$ Then by the `pigeonhole principle' there is some string $\nu$ (possibly empty) such that $|\nu\tau_{-}| = n$ and $\nu\tau_{-}$ occurs infinitely many times. But the formation rule for $EM_n(Z)$ then yields infinitely many occurrences of both $\tau_{-}0$ and $\tau_{-}1$. In particular, $\tau$ itself occurs infinitely many times, contrary to assumption. Thus it must be the case that
 $|\tau_{-}| \ge n+1. \hfill\square$

\begin{proposition} Let X be a DeBruijn string of order $n+1$ with initial n-string $\sigma$, say X = $\sigma{Y}\sigma.$ Then $EM_n(X) =
Y\sigma{Y}\sigma\dots.$
\end{proposition}

This is an easy consequence of Proposition 3.1.

\begin{theorem} For any key Z, $EM_n(Z)$ will eventually produce a DeBruijn string of order $n+1.$ After that, the sequence will repeat as if the DeBruijn string had been the key.
\end{theorem}

Proof: Let $\tau$ be a string of length $n + 1$. We shall say that a natural number $p$ is a {\it period} of $\tau$ if there are infinitely many non-overlapping substrings X of $EM_n(Z)$ with $|X| = p$, and that begin and end with $\tau$ with no other $\tau$ in between. Since $EM_n(Z)$ is ultimately periodic, we may choose and fix a string $\tau$ having a maximum possible period. Let W denote the (infinite) tail of $EM_n(Z)$ that begins at the earlier of the two occurrences of $\tau$ in one such X.

Let $\sigma$ be a string of length $n+1$ occurring in W such that the substring extending from the first occurrence of $\sigma$ to the second is as short as possible. If more than one string contends for this distinction, choose the earliest occurring one in W.

We shall argue now that $\sigma = \tau$, in particular, $\sigma$ is the initial string of W. Note that there must be a distinct instance of $\sigma_{-}$ between $\sigma$ and its next occurrence. In other words W must contain
$$
\sigma \dots x\sigma_{-} \dots \sigma
$$
for some binary digit $x$. (As in \cite{me}, we understand in this and similar diagrams that the indicated substrings are distinct and occur in the indicated order, but that overlap may occur.) If the predecessors of the two $\sigma$s were the same, then the first one would have to be the initial string of W, since $\sigma$ was chosen to occur as early as possible in W. Thus we may assume the two $\sigma$s have opposite predecessors. In that case, one of the two predecessors must be $x$. If it were the predecessor of the second $\sigma$ then the substring $x\sigma_{-} \dots x\sigma_{-}$ would be shorter than $\sigma \dots \sigma$, contrary to assumption. Thus $x$ must match the predecessor of the first $\sigma$. But in that case we reach a similar contradiction unless the first $\sigma$ were the initial string of W.

By the forgoing, every string of length $n+1$ that occurs in W must recur at equally spaced intervals, and that spacing is the same as for $\sigma$, the initial string of W. Thus $W = \sigma_{-}Y\sigma_{-}Y\dots$ for some string Y. Let $X = \sigma_{-}Y\sigma_{-}.$ Let $\mu$ be any given string of length $n+1.$ By Proposition 3.2, some instance of $\mu$ must begin in $\sigma_{-}Y.$ Thus, $\mu$ is a substring of X. It follows that X contains every given string of length $n+1$ exactly once, and it is therefore a DeBruijn string.$\hfill\square$
\vskip .1 in
The previous result is impractical for generating DeBruijn sequences from an arbitrary seed since it gives no upper bound on the number of steps required. On the other hand, it is easy to produce seed strings $Z$ for which $EM_n(Z)$ is immediately periodic.

\begin{theorem} Let Z be any string ending in $n$ zeros that contains every n-string at least once. Then the initial $2^{n+1}$-string of $EM_n(Z),$ together with the $n$ zeros at the end of Z, is a DeBruijn string of order $n+1.$
\end{theorem}

Proof: \cite{fred} presents an algorithm that produces an Eulerian circuit of ${\cal B}^n$ from a directed tree in ${\cal B}^n$ rooted at $\bar 0_n$.(The corresponding Hamilitonian circuit of ${\cal B}^{n+1}$ is then a DeBruijn string of order $n+1.$) Consider the subgraph, $T$, of ${\cal B}^n$ formed by connecting the vertex of an n-string $\sigma$ to the vertex of the n-string that follows the {\it last} occurrence of $\sigma$ in Z. It is obvious that every n-string can be connected to $\bar 0_n$: just follow along Z and jump to the last occurrence whenever necessary. Thus $T$ is a tree rooted at $\bar 0_n$. Now apply ALGORITHM 1 on page 200 of \cite{fred}. (One must replace the $n$ in \cite{fred} with $n+1.$) It is only necessary to note that the steps of the algorithm are identical with those of $EM_n(Z).$ $\hfill\square$
\vskip .1 in
With a slight modification, the algorithm can be adapted to handle arbitrary seeds: Start with $\bar 0_n$. Every time an n-string tail is encountered that does not appear in the seed, pre-pend it to the beginning of the seed and continue as if the new longer seed had been the one given. The algorithm is very fast, and is simple enough to use by hand.  For example, with $n = 3$ and an empty seed, it produces the order 4 DeBruijn string 0001111010110010000. (In the course of the algorithm, the original empty seed turns into 111011001000.)
\vskip .1 in
We return now to the study of the usual EM sequence (i.e., with $h$ equal to the identity function,) but with a nonempty seed string Z. The resulting sequences EM(Z) share some of the important features of the usual EM sequence. For example, we have

\begin{proposition} Let $m_0$ be the length of the longest string that recurs in Z. Let $\sigma$ be the match string of EM(Z) at time $n.$ If $|\sigma| > m_0$ and the match length $m_n$ reaches a new record value at time $n+1$, i.e., $m_{n+1} > m_k, k = 1,2, \dots, n$, then $\sigma$ is a head of ZEM(Z).
\end{proposition}

The proof is essentially the same as that of Proposition 4.1 of \cite{me}.

Let $T_n$ be the first time the match length reaches $n$. Then we have
\begin{equation}
T_n \le 2^n + n - |Z| < 2^{n+1} - |Z|, n > m_0.
\end{equation}
To see this, note that if $T_n > k$ then the finite string $ZEM(Z)_k$ contains no repeated strings of length $n$. Thus $k + |Z| - n + 1 \le 2^n$ and the desired inequality follows by taking $k = T_n - 1.$

Like EM, the strings $EM(Z)$ are always transitive. Let $C_n$ be the {\it cover time} of strings of length $n$, i.e., the smallest k for which every string of length $n$ is a substring of $EM(Z)_k.$ We shall show $EM(Z)$ is transitive by obtaining a crude upper bound for the $C_n$. For a given positive constant $C$, define a function $f$ on the natural numbers inductively by setting $f(1) = C,$ and $f(n+1) = 2^{f(n)+|Z|+2}.$

\begin{proposition} The constant $C$ can be chosen depending only on Z such that $C_n \le f(n).$
\end{proposition}

Proof: We proceed by induction on $n$. It is easy to see that EM(Z) cannot consist entirely of ones or entirely of zeros. Thus we can handle the case $n = 1$ merely by taking the constant $C$ sufficiently large.

Put $N = f(n)$ and assume as an inductive hypothesis that every string $\tau$ of length $n$ occurs at least once in $EM(Z)_N.$ Fix any such $\tau$ and let $\sigma$ be the prefix of ZEM(Z) that ends at one of the occurrences of $\tau$ in $EM(Z)_N.$ By Proposition 3.4, both $\sigma0$ and $\sigma1$ occur in the initial segment of ZEM(Z) having length
$T_{|\sigma|+1}.$ Hence, both $\tau0$ and $\tau1$ occur in the initial segment of EM(Z) having length $T_{|\sigma|+1} - |Z|.$ But $|\sigma| \le N + |Z|,$ so by (3.2) $T_{|\sigma| + 1 } \le 2^{|\sigma|+ 2} \le 2^{N + |Z| + 2} = f(n+1).$

Since $\tau$ was arbitrary, every string of length $n+1$ occurs at least once in the initial segment of EM(Z) of length $f(n+1),$ i.e., $C_{n+1} \le f(n+1).$ This completes the inductive step, and the proof.

Next, we show that Z is uniquely determined by EM(Z).

\begin{theorem} Let $Z_1$ and $Z_2$ be binary strings. Then if $EM(Z_1) = Y = EM(Z_2)$, the strings $Z_1$ and $Z_2$ are equal.
\end{theorem}

Proof: Assume $Z_1 \ne Z_2.$ We may assume that $|Z_1| \ge |Z_2|$ and that
 $Z_1$ is not the empty string. Let $X_1 = Z_1EM(Z_1)$ and $X_2 = Z_2EM(Z_2)$.

 Let $\sigma$ be the head of $Y$ of length $k$. Since $Y$ is not periodic, the distance $d(k)$ between the first two occurrences of $\sigma$ in $Y$ tends to infinity as $k$ tends to infinity. Fix $k$ with $d(k) \ge |Z_1| + 1.$

 The first matches of $\ _-Z_1\sigma$ and $Z_1\sigma$ in $X_1$ must follow one of the following two patterns:
 \begin{enumerate}
 \item[i]
 $$
 Z_1{\sigma}x \ \ \dots \ \ \ _-Z_1{\sigma}x^{\prime} \ \ \dots \ \ Z_1{\sigma}x^{\prime}
 $$
 \item[ii]
 $$
  Z_1{\sigma}xy \ \ \dots \ \ Z_1{\sigma}x^{\prime} \ \ \dots \ \ \ _-Z_1{\sigma}xy^{\prime}
 $$
 \end{enumerate}

See, e.g., Theorem 4.6 of \cite{me} or Lemma 2 of \cite{sutner}. One may check that these results generalize to the seeded case provided all match lengths involved in the arguments exceed the length of the seed. (In the language of \cite{me},  an excursion begins in case (ii) at the middle sequence in the display, and ends at the last sequence. The key point here is that the match of $\ _-Z_1{\sigma}$ occurs before the next appearance of $Z_1\sigma.$)

 In case (i), the choice of $k$ ensures that the string denoted $\ _-Z_1{\sigma}x^{\prime}\dots Z_1{\sigma}x^{\prime}$ lies entirely inside $Y$, i.e., does not overlap with the seed. Thus, in $X_2$ we have
 $$
 Z_2{\sigma}x \ \ \dots \ \ \ _-Z_1{\sigma}x^{\prime} \ \ \dots \ \ Z_1{\sigma}
 $$
 Since $|Z_1| \ge |Z_2|$ and $Z_1$ and $Z_2$ are not equal, there can be no other occurrence of $Z_1\sigma$ in the range shown. But then the match of $\ _-Z_1\sigma$ at the end would produce a next digit of $x$ rather than $x^{\prime}$, a contradiction.

 In case (ii), we replace $\sigma$ with the next longer initial string ${\sigma}x$. Then at the next occurrence of $Z_1{\sigma}x$ in $X_1$ we are back in case (i).
 $\hfill\square$

\section{Double Helices}

A given DeBruijn string $X$ of order $n$ induces a Hamiltonian circuit of the DeBruijn graph ${\cal B}^n$ that begins at the node corresponding to the initial n-string of $X$. For example, the order 3 string 0011101000 induces the circuit $001 \to 011 \to 111 \to 110 \to 101 \to 010 \to 100 \to 000 \to 001$ of ${\cal B}^3.$ If the edges of such a circuit are removed from
${\cal B}^n$, the remaining graph contains at least 3 connected components, including the loops from the all zero n-string back to itself, and from the all 1 n-string back to itself. For the order 3 example just given, there is a single additional component consisting of the loop $001 \to 010 \to 101 \to 011 \to 110 \to 100 \to 001.$

In cases where there are exactly 3 connected components (i.e., the minimum number possible,) we shall call the original DeBruijn string $X$ a {\it double helix.} The long cycle in the graph formed by removing the edges of the path of $X$ from the DeBruijn graph will be called the {\it message loop of X,} and the corresponding string the {\it message of X.} We also call a DeBruijn cycle a double helix if the associated DeBruijn string is a double helix.

The message of an order $n$ double helix is not itself a DeBruijn string since it does not include the all-zero and all-one strings, $\bar{0}_n$ and $\bar{1}_n$. It does contain every other binary string of length $n$. Such a string can be converted to a {\it bona fide } DeBruijn string by inserting an extra 0 at the position of $\bar{0}_{n-1}$ and an extra 1 at the position of $\bar{1}_{n-1}.$ We shall call a binary string that contains every string of length $n$ exactly once, except for one or both of the all-zero and all-one strings, a {\it depleted DeBruijn string.}

Every DeBruijn string of order 3 or less is a double helix, but this is not true for higher orders. For example, 1111000010011010111 is a double helix of order 4, but the DeBruijn string 1111000011010010111 is not a double helix.

If a double helix $Z$ of order $n$ is used as a seed for either $EM$ or $EM_n$, then it is nearly correct to say that the seeded sequence begins by extracting the message of Z. More precisely, we have:

\begin{theorem}
Let $Z$ be an order $n$ double helix whose initial $n-$string $\tau$ is neither the all-zero nor the all-one string. Let $Y$ be the initial $2^n - 1$ string of either $EM_n(Z)$ or $EM(Z).$ Then the path of $\tau_-Y$ is the message loop of $Z$, and $ZY$ is a depleted DeBruijn string of order $n + 1.$
\end{theorem}

Proof: The theorem is vacuous if $n = 1$ and is readily checked for $n = 2$, so we may assume that $n \ge 3.$ Let ${\sigma}xy$, where $x$ and $y$ are binary digits, be the initial $n+1$ string of $Z$. Then $\sigma$ is also the terminal $n-1$ string of $Z$, and $xy^{\prime}$ are the first two digits of $Y$. Now $\tau = {\sigma}x$ lies on the message loop $L$ of $Z$ in $\mathcal{B}^n$, and
${\sigma}x \to {}_-{\sigma}xy^{\prime}$ is an edge of this loop. By induction, each successive edge of the path of $ZY$ is an edge of $L$, until the entire loop shall have been traversed. Completion of $L$ starting from ${\sigma}x$ with a return to ${\sigma}x$ on the last step requires $2^n - 1$ steps, since the path visits every $n-$string except $\bar0_n$ and $\bar1_n$ exactly once.

Every $n+1$ string except $\bar0_{n+1}, \bar1_{n+1}$ occurs exactly once in $ZY$. Therefore, $ZY$ is a depleted DeBruijn string of order $n+1$.$\hfill\square$

The method of Theorem 4.1 produces depleted DeBruijn strings of a special type: their first half is a DeBruijn string of the next lower order. They are not themselves depleted double helices, in general.

Double helices do not appear to be particularly rare among DeBruijn strings. For example, 4 of the 16 order 4 DeBruijn cycles are double helices, and 840 of the 2048 order 5 DeBruijn cycles are double helices. We shall show there are double helices of all orders by showing that all DeBruijn strings generated by means of linear recurrences are double helices, but this only accounts for 5 of the 840 double helices of order 5. It would be of interest to find other algorithms that reliably generate double helices of all orders, as well as to count the exact number occurring at each order.

Consider a linear recurrence of the form
\begin{equation}
a_k = c_1 a_{k-1} + c_2 a_{k-2} + \dots + c_{n}a_{k-n},\ \ k = 0, 1, \dots,
\end{equation}
where all operations and elements are those of $\mathbb{Z}_2$, the field of integers modulo 2. If a set of $n$ coefficients $c_i$ is given, as well as a choice of initial values $a_{-1}, a_{-2}, \dots a_{-n}$, then the recurrence determines a unique infinite sequence of binary digits $a_0, a_1, \dots$ Since there are $2^n$ possible $n-$tuples of binary digits, the sequence must eventually become periodic with period at most $2^n.$ The maximal period is, in fact, $2^n - 1$, since a sequence of n zero digits can only be followed by zero digits.

Binary sequences produced by a linear recurrence that have maximal period are called PN sequences in \cite{golomb}. A full period of a PN sequence contains all binary sequences of length $n$ except for the all zero sequence. Inserting an extra zero at the beginning of the unique sequence of $n - 1$ zeros in such a cycle produces a DeBruijn cycle. We shall call these cycles and their associated sequences {\it linear DeBruijn } cycles and sequences.

Clearly a necessary condition for a PN sequence is that coefficient $c_n$ should be non-zero. It is convenient to introduce an additional coefficient $c_0 = 1,$ and to define the {\it characteristic polynomial}
$$
f(x) = \sum_{j=0}^n c_j x^j.
$$
It is shown in \cite{golomb} that characteristic polynomials of PN sequences are necessarily irreducible in the polynomial ring $\mathbb{Z}_2[x]$, and that the period $p$ of any sequence generated by a linear recurrence, for which the characteristic polynomial is irreducible, is determined as the smallest value of $p$ such that the polynomial $1 - x^p$ is divisible by the characteristic polynomial. Furthermore, PN sequences exist for every order $n.$ It therefore follows from our next result that double helices exist for all orders $n$.

\begin{theorem}Every linear DeBruijn sequence is a double helix.
\end{theorem}

Proof: Let $a_n$ be a PN sequence determined by the recurrence (4.1) and having characteristic function $f(x)$. Let $b_n$ be the sequence determined by the linear recurrence
\begin{equation}
b_k = c_1 b_{k-1} + c_2 b_{k-2} + \dots + c_{n}b_{k-n} + 1,\ \ k = 0, 1, \dots,
\end{equation}
with initial condition $b_{-1} = 1, b_{-2} = b_{-3} = \dots = b_{-n} = 0.$ The paths in $\mathcal{B}^n$ induced by the $a_k$ sequence and the $b_k$ sequence leave each vertex by the opposite edge, so it suffices to show that the period of the latter sequence is $q = 2^n - 1.$

Let $G(x)$ be the generating function of the $b_k$ sequence, i.e.,
$$
G(x) = \sum_{k=0}^{\infty}b_k x^k.
$$
Substituting for $b_k$ the expression on the right side of (4.2) and interchanging orders of summation, we obtain
$$
G(x) = \frac1{1 - x} + \sum_{i=1}^n c_i x^i(b_{-i}x^{-i} + \dots + b_{-1}x^{-1} + G(x)),
$$
whence
$$
G(x)f(x) = \frac1{1 - x} + \frac{1 - f(x)}{x} = \frac{1 + xf(x) + f(x)}{x(1-x)}.
$$
(Recall that operations are done in $\mathbb{Z}_2[x].$)
Now since $f(0) = 1,$ the last written numerator is divisible by $x$, i.e.
\begin{equation}
h(x) = \frac{ 1 + xf(x) + f(x)}{x} = f(x) + \sum_{i=1}^nc_ix^{i-1}
\end{equation}
is a polynomial that also satisfies
$$
G(x)(1-x)f(x) = h(x).
$$
On the other hand, since $b_k$ has period $q$, we have
$$
G(x) = \frac{g(x)}{1 - x^q}
$$
for some polynomial $g(x)$. Combining these results we have $g(x)(1-x)f(x) = h(x)(1 - x^q).$ Since $f(x)$ is irreducible, it must divide either $h(x)$ or $ 1 - x^q.$ But by (4.3) the first possibility cannot occur, since then $f$ would divide a polynomial of lower degree. It follows that $f$ divides $1 - x^q$, and hence $q \ge 2^n - 1.$ But since the $b_n$ sequence never visits the all one n-string, its period is exactly $2^n - 1$, as required. (Since a PN sequence does visit $\bar 1_n$, it is necessary that $c_1 + c_2 + \dots + c_n = 0$ to avoid getting trapped there. Thus the $b_n$ sequence cannot visit $\bar 1_n$ unless it starts there.) The path induced in ${\cal B}^n$ by the $b_n$ sequence takes the form
$\bar 0_{n-1}1 \to \dots \to 1\bar 0_{n-1} \to \bar 0_n \to \bar 0_{n-1}1.$ ($2^n -1$ steps.) If the original PN sequence is converted to a DeBruijn sequence by inserting a zero, then in the message loop the last two steps are replaced by the single step $1\bar 0_{n-1} \to \bar 0_{n-1}1.$ Thus the message loop omits both $\bar 0_n$ and $\bar 1_n$. $\hfill\square$

\end{document}